\documentclass[12pt]{article}
\usepackage{amsmath,amsfonts}
\begin{document}

\baselineskip 20pt
\date{\ April 30, 2014}
\title{Weak-star convergence and a polynomial approximation problem}

\author{Arthur A.~Danielyan}

\maketitle

\begin{abstract}
\noindent Let $E$ be an arbitrary subset of the unit circle $T$ and
let $f $ be a function defined on $E$. When there exist polynomials
$P_n$ which are uniformly bounded by a number $M>0$ on $T$ and
converge (pointwise) to $f$ at each point of $E$? We give a
necessary and sufficient description of such functions $f$.
The necessity part of our result, in fact, is a classical theorem of
S.Ya. Khavinson, while the proof of sufficiency uses the method that
has been recently applied in particular in the author's solution of an approximation
problem proposed by L. Zalcman.

\end{abstract}

\begin{section}{Introduction and the result.}

For a compact subset $K$ of $\mathbb C$ we denote by $C(K)$ the
space of all continuous complex valued functions on $K$. Denote by
$\Delta$ and $T$ the open unit disc and the unit circle, respectively.

Let $F$ be a closed subset of the unit circle $T$ and
let $f \in C(F)$. {\it The problem of uniform
approximation of $f$ on $F$ by polynomials $P_n$ that are
uniformly bounded on $\Delta$} has been investigated in \cite{Zalc} and  \cite{Dan}. In case when $F$ is a closed arc of $T$, it has been
solved in \cite{Zalc}. In the general case the problem has been settled by the following theorem (see \cite{Dan}).

\vspace{0.25 cm}

{\bf Theorem A.}\ {\it Let $F$ be a closed subset of $T$. In order
that a function $f \in C(F)$ be uniformly approximable on $F$ by
polynomials $P_n$ such that $|P_n(z)| \leq M, \ z \in \Delta, \ n
= 1,2,...$, it is necessary and sufficient that the following
conditions be satisfied:

(a) $|f(z)|\leq M$, $z \in F$, and

(b) There exists a function $g$ analytic on $\Delta,  \ |g(z)|
\leq M, \ z \in \Delta $, such that $$f(e^{i \theta})= \lim_{r \to
1} g(re^{i\theta})$$ for almost all $e^{i \theta} \in F$.}

\vspace{0.25 cm}

The paper  \cite{Zalc} suggests also further research ideas and
formulations of open approximation problems.
 For example,  the above approximation problem
 can be modified {\it ``by altering variously the set on which one approximates, the sense in which approximation is required to hold, and the precise conditions of boundedness"}     (see \cite{Zalc}, p. 380).
In the present paper we address
the case when the condition of uniformly boundedness is the same as in
the original problem, while the set on which one approximates is an
arbitrary subset of $T$ and the approximation is merely pointwise.
The precise formulation of the
problem is the following.

\vspace{0.25 cm}

{\bf Problem 1.} {\it Let $E$ be an arbitrary proper subset of $T$ and
let $f $ be a function defined on $E$. Find necessary and sufficient conditions on $f$ in order there
exist polynomials
$P_n$ which are uniformly bounded by $M>0$ on $T$ and
converge (pointwise)  to $f$ at each point of $E.$}

\vspace{0.25 cm}

The (uniformly) boundedness condition in Problem 1 is important, because if it is dropped,
the appropriate question will have the following trivial answer:
The function $f$ on $E$ is a limit of a sequence of polynomials $P_n$ if and only if there exists a sequence of
continuous functions on $T$ which converges (pointwise) to $f$ at each point of $E$.

Note that Theorem 2 of the paper  \cite{Zalc} solves another typical bounded approximation
problem with an $L^p$ norm  boundedness  condition (on $T$) for the polynomials which approximate the function $f$ {\it pointwise and a. e.  on an arbitrary (measurable) subset} $E$  of $T$. The precise formulation of this result is the following.

\vspace{0.25 cm}

{\bf Theorem B.} {\it Let $E \subset T$ be a set of positive measure and let $q \geq 1$. A function $f$ on $E$ is the (pointwise, almost everywhere) limit of polynomials $P_n$ satisfying $||P_n||_q \leq M$ if and only if there exists a function $g \in H^q$, $||g||_q \leq M$, such that
$f(z)=g(z) \equiv  \lim_{r \to
1} g(rz)$ for almost all $z \in E.$}

\vspace{0.25 cm}

Let the functions $f_n, \ n= 1, 2, ...,$ and $g$ be elements of  $L^\infty$ on $T$. Recall that the sequence $\{f_n\}$  converges weak-star to $g \in L^\infty$
if

\begin{equation}
\lim_{n\rightarrow \infty}\int_T f_n(z)Q(z)dz = \int_{T}
g(z)Q(z)dz
\end{equation}
for any function $Q \in L^1$ on $T.$

By Fatou's theorem, to each bounded analytic in $\Delta$ function corresponds its (radial) boundary function defined almost everywhere on $T$.

The following classical theorem is due to S. Ya. Khavinson (see \cite{Priv}, p. 171).

\vspace{0.25 cm}

{\bf Theorem C.}\ {\it Let the sequence  $\{f_n\}$  of analytic in $\Delta$ functions be uniformly bounded.
If  $\{f_n\}$ converges to a limit function $g$ in $\Delta$,
then the boundary functions of
$\{f_n\}$ converge weak-star to the boundary function of
$g$.}

\vspace{0.25 cm}

Now we formulate the solution of Problem 1, which is the aim of this paper.

\vspace{0.25 cm}

{\bf Theorem 1.}\ {\it Let $E$ be an arbitrary subset of $T$ and let  $f$ be a function on $E.$
There exists a sequence of  polynomials $\{P_n\}$ converging (pointwise) to $f$ at each point of $E$ and satisfying the condition $|P_n(z)| \leq M, \ z \in \Delta, \ n
= 1, 2,...$,
if and only if there exist a sequence of continuous functions  $\{f_n\}$ on $T$ and a bounded analytic function $g$ in  $\Delta$, such that:

(i) $|f_n(z)| \leq M, \ z \in T, \ n
= 1, 2,...$, and $\{f_n\}$ converges (pointwise) to $f$ at each point of $E$;  and

(ii) $\{f_n\}$ converges weak-star to the boundary function of $g$.}

\vspace{0.25 cm}

The following corollary describes a particular class of approximable functions on an arbitrary measurable subset $E$ of $T$.

\vspace{0.25 cm}

{\bf Corollary 1.}\ {\it Let $E$ be an arbitrary measurable subset of $T$ and let  $f$ be a function on $E$.
Assume that there exists
a bounded analytic function $g$ in  $\Delta$, the boundary function of which coincides with $f$ a. e. on $E$.
Assume furthermore that  there exist a sequence of continuous functions  $\{f_n\}$ on $T$ such that:

(i) $|f_n(z)| \leq M, \ z \in T, \ n
= 1, 2,...$, and $\{f_n\}$ converges (pointwise) to $f$ at each point of $E$;  and

(iii) $\{f_n\}$ converges a. e. on $T \setminus E$ to the boundary function of $g$.

\noindent Then there exists a sequence of  polynomials $\{P_n\}$ converging (pointwise) to $f$ at each point of $E$ and satisfying the condition $|P_n(z)| \leq M, \ z \in \Delta, \ n
= 1, 2,...$.}

\vspace{0.25 cm}

To prove Corollary 1, first note that the condition {\it (i)} is the same in both Theorem 1 and Corollary 1. If  the conditions of Corollary 1 are satisfied, then obviously, the sequence  $\{f_n\}$  converges to the boundary function of $g$ a. e. on $T$, and by the Lebesgue dominated convergence theorem, (1) is satisfied (for any $Q \in L^1$ on $T$). Thus, the condition {\it (ii)} of Theorem 1 is satisfied as well, and Corollary 1 follows from Theorem 1.

The proof of Theorem 1 follows the approach applied in particular in \cite{Kol}  and    \cite{Dan}; for further details  and recent applications  see \cite{Dan} and  \cite{Dan2}.

\end{section}

\begin{section}{Proof of Theorem 1}

The necessity of the condition {\it (i)} is obvious.
Since the sequence of polynomials $\{P_n\}$ is uniformly bounded on $\Delta$, without loss of generality we may assume
it
converges uniformly on compact subsets of $\Delta$ to a bounded analytic function $g.$
By the Theorem C the sequence $\{P_n\}$ converges weak-star to the
boundary function of $g.$ The necessity of {\it (ii)} is proved.

Now we present the proof of the sufficiency of the conditions of
Theorem 1.
Let $C=C(T)$ and let $A$ be the disk algebra. By the Riesz
theorem, a continuous linear functional
$\Psi$ on the quotient space $C/A$ can be represented as
$$\Psi(h+ A)= \int h d\mu,$$
where $h+A = \{h+u: u \in A\}$ is an element of the space $C/A$
and $\mu$ is a complex Borel measure on $T$ orthogonal to $A$ (cf.
\cite{Koos}, p. 192
).  By the F. and M. Riesz theorem
$d\mu=G(z)dz$ for some $G \in H^1,$ and, thus

$$\Psi(h+ A)= \int h d\mu = \int_{T} h(z)G(z)dz.$$

Since $\{f_n\}$ converges weak-star to the boundary function of $g,$
one can apply the relation (1) for the function $G \in H^1 \subset L^1,$
which (combined with $gG \in H^1$) implies

$$\lim_{n \rightarrow \infty} \Psi (f_n +A)=  \lim_{n\rightarrow \infty}\int_T f_n(z)G(z)dz = \int_{T}
g(z)G(z)dz= 0.$$

Hence, $\{f_n + A\}_{n=1}^{\infty}$ converges weakly to zero in
$C/A.$ Obviously the same is true for the sequence $\{f_n + A\}_{n=m}^{\infty}$ for
 each fixed natural number $m.$ By a theorem of Mazur (see \cite{Yos}, p. 120), there exists a sequence of finite convex
linear combinations of elements of $\{f_n + A\}_{n=m}^{\infty}$  which converges
to zero in the norm of $C/A.$ Thus for
each natural number $m$ there exist numbers $\alpha_n^{(m)} \geq 0
\ \ (n = 1, 2, ..., k_m)$, such that $
\sum_{n=m}^{k_m}\alpha_n^{(m)} =1$ and the quotient norm of the
linear combination $\sum_{n=m}^{k_m}\alpha_n^{(m)}(f_n + A)=
\sum_{n=m}^{k_m}\alpha_n^{(m)} f_n +A$ is less than $1/m$. By the
definition of the quotient norm there exists a
function $u_m \in A$ such that
\begin{equation}
|\sum_{n=m}^{k_m}\alpha_n^{(m)} f_n(z) - u_m(z)|< 1/m, \ \  z \in
T.
\end{equation}

In particular
\begin{equation}
|u_m(z)|< |\sum_{n=m}^{k_m}\alpha_n^{(m)} f_n(z)|+ 1/m \leq
\sum_{n=m}^{k_m}\alpha_n^{(m)} M + 1/m = M +1/m  , \ \ z \in T.
\end{equation}

It is easy to see that $\{u_n\}$ converges to $f$ at each point of $E.$

Let $Q_m$ be a polynomial such that $|u_m(z) - Q_m(z)|<1/m, \ z
\in T$. By (3) we have $|Q_m(z)|< M+ 2/m, \ z \in T$. The
polynomials $$P_m(z)= \frac{M}{M+ 2/m}\ Q_m(z) $$ are uniformly
bounded on $T$ by $M$ and converge to $f$ on $E.$
The proof is over.

\end{section}

\begin{minipage}[t]{6.5cm}
Arthur A. Danielyan\\
Department of Mathematics and Statistics\\
University of South Florida\\
Tampa, Florida 33620\\
USA\\
{\small e-mail: adaniely@usf.edu}
\end{minipage}

\end{document}